\documentstyle{amsart}
\input{bull-art}
\newcommand{\ha}{\frac{1}{2}}
\newcommand{\ol}{\overline}
\newtheorem{theo}{Theorem}[section]

\newtheorem{lem}[theo]{Lemma}

\begin{document}
\def\currentvolume{26}
\def\currentissue{2}
\def\currentyear{1992}
\def\currentmonth{April}
\def\copyrightyear{1992}
\def\currentpages{322-328}
\title[Sampling and interpolation in the Bargmann-Fock 
spaces]
{ Density theorems for sampling and
interpolation \\ in the Bargmann-Fock space}
\author{Kristian Seip}
\date{August 9, 1991 and, in revised form, November 19, 
1991}
\address{Division of Mathematical Sciences, Norwegian 
Institute
of Technology, N-7034 Trondheim-NTH, Norway}
\email{seip@@imf.unit.no}
\maketitle
\begin{abstract}
We give a complete description of sampling and 
interpolation in
the Bargmann-Fock space, based on a density concept of 
Beurling. 
Roughly speaking, a discrete set is a set of sampling
if and only if its density in every part of the plane is 
strictly
larger than that of the von Neumann lattice, and similarly, 
a discrete set is a set of interpolation if and only if
its density in every part of the plane is strictly smaller 
than
that of the von Neumann lattice. 
\end{abstract}
\section{Introduction and results}
The work presented in this announcement is based on 
Beurling's 
lectures on balayage of Fourier-Stieltjes transforms and
interpolation for an interval on $\bold R^1$ \cite{be}. 
We observe that Beurling's problems concerning functions 
of exponential
type have natural counterparts for
functions of order two, finite type and find that, indeed, 
so have his
main results. The most interesting part, however, is that 
Beurling's ideas
are applicable also in the  Hilbert space setting, yielding
a complete  description of sampling and interpolation in the
Bargmann-Fock space. The simplicity of
these results is quite remarkable when compared to the 
situation
in the Paley-Wiener space 
(the corresponding Hilbert space of functions of 
exponential type)
and to the extensive literature on nonharmonic Fourier 
series and,
in particular, Riesz bases of complex exponentials [20].
  
This research is motivated by a  
recent development in signal analysis and applied 
mathematics,
which was initiated by Daubechies, Grossmann, and Meyer 
\cite{dagr,da2,dagrme}.
Their work inspired us to search for a general 
characterization of the 
information needed to represent signals, as functions in 
the 
Bargmann-Fock space. Our results can be seen as sharp 
statements about
the Nyquist density and its meaning in this context.

In order to describe more precisely the problems to be 
considered,
a few definitions are needed. For  $\alpha>0$, 
let $d\mu_{\alpha}(z)=(\alpha/\pi)e^{-\alpha
|z|^2}dxdy$, $z=x+iy$, and define the Bargmann-Fock space 
$F_{\alpha}^2$
to be the collection of entire functions $f(z)$ for which
\[
\| f \|_{2} = \| f \|_{\alpha,2} =
\int_{{\bold C}} |f(z)|^2 d\mu_{\alpha}(z)<\infty.
\]
$F_{\alpha}^{2}$ is a Hilbert space with reproducing kernel
$K(z,\zeta)=e^{\alpha \ol{z}\zeta}$;
i.e., for every $f\in F_{\alpha}^{2}$ we have
\[
f(z)=\langle f,K(z,\cdot)\rangle =
\int_{{\bold C}}f(\zeta)\ol{K(z,\zeta)}\ 
d\mu_{\alpha}(\zeta).
\label{repr}
\]
The normalized reproducing kernels, 
$k_{\zeta}(z)=K(\zeta,\zeta)^{-1/2}K(\zeta,z)$, can be 
view-\linebreak ed as the
natural (well-localized) building blocks of $F_{\alpha}^2$.
They correspond, via the Bargmann transform,
to the canonical coherent states of quantum mechanics
and to Gabor wavelets in signal analysis. This relation is 
the reason
for the importance of the Bargmann-Fock space; see 
\cite{fo} for general information and \cite{dagr} for more 
background on the problems  treated here. 

We say  that a discrete set $\Gamma$ of complex numbers is 
a 
{\em set of sampling} for $F_{\alpha}^{2}$ if there exist
positive numbers $A$ and $B$ such that
\begin{equation}
A \| f \|_2^2 \leq  \sum_{z\in \Gamma}e^{-{\alpha} 
|z|^2}|f(z)|^2
\leq B \| f \|_2^2
\label{frame}
\end{equation}
for all $f\in F_{\alpha}^{2}$. If to every $l^2$-sequence
$\{a_{j}\}$ of complex numbers there exists an $f\in
F_{\alpha}^{2}$ such that 
$e^{-\alpha|z_{j}|^2/2}f(z_{j})=a_{j}$
for all $j$, then $\Gamma=\{z_{j}\}$ is said to be a {\em 
set of interpolation}
for $F_{\alpha}^{2}$. A set of sampling corresponds, 
in the terminology of \cite{dusc}, to a {\em frame} of 
coherent states. A set
of both sampling and interpolation (which does not exist) 
would correspond to a {\em Riesz basis}
of coherent states.

With a view to applications in physics and signal analysis,
Daubechies and Grossmann posed  the problem of 
finding the lattices $z_{mn}=ma+inb$, $m,n\in {\bold Z}$, 
that are sets
of sampling \cite{dagr}. 
They proved that a lattice could be a set of sampling only
if $ab<\pi/\alpha$ and conjectured  this condition also to 
be 
sufficient. For $ab=\pi/(\alpha N)$, $N$ an integer 
$\geq 2$, they found (\ref{frame}) to hold by providing 
explicit 
expressions for the optimal constants $A,B$. Daubechies 
was later able to 
show that
a lattice is a set of sampling whenever $N^{-1}<0.996$ 
\cite{da2}.

We prove that the density criterion of the 
Daubechies-Grossmann 
conjecture applies not only to lattices, but to {\em 
arbitrary} discrete 
sets. We should add here that the conjecture was proved 
independently by
Lyubarskii \cite{ly} and by Wallst\'{e}n and the author 
[19].

For the description to be given of sets of sampling and 
interpolation, 
we need Beurling's density concept as
generalized by Landau \cite{la}.
We consider then {\em uniformly discrete sets}, i.e., 
discrete sets 
$\Gamma=\{z_{j}\}$ for which $q=\inf_{j\neq k} 
|z_{j}-z_{k}|\break
>0$. 
We fix a compact set $I$ of
measure $1$ in the complex plane, whose boundary has 
measure 0.
Let $n^{-}(r)$ and $n^{+}(r)$ 
denote, respectively, the smallest and largest number of 
points 
from $\Gamma$ to be found in a translate of $rI$.
We define the lower and upper uniform densities of $\Gamma$
to be
\[
D^{-}(\Gamma)=\liminf_{r\rightarrow\infty}%
\frac{n^{-}(r)}{r^2}
\quad\mbox{and}\quad
D^{+}(\Gamma)=\limsup_{r\rightarrow\infty}\frac{n^{+
}(r)}{r^2},
\]
respectively. It was proved by Landau that these limits are
independent of $I$.

Our main theorems are the following (the sufficiency
part of the theorems were obtained in collaboration with 
Wallst\'{e}n [19].

\begin{theo}
A discrete set $\Gamma$ is a set of sampling for
$F_{\alpha}^{2}$ if and only if it can be expressed as a 
finite
union of uniformly discrete sets
and contains a uniformly discrete
subset ${\Gamma}'$ for which $D^{-}({\Gamma}')>\alpha/\pi$.
\label{sampling}
\end{theo}

\begin{theo}
A discrete set $\Gamma$ is a set of interpolation for
$F_{\alpha}^{2}$ if and only if it is uniformly discrete
and $D^{+}({\Gamma})<\alpha/\pi$.
\label{interpol}
\end{theo}

{\em Remark} 1.  
Decomposition and interpolation theorems for general 
discrete sets were
obtained in \cite{japero}, however, without any indication 
of a 
critical density. The results in \cite{japero} appear as 
part of
a certain trend in harmonic
analysis, and the analogy to the theory of nonharmonic
Fourier series does not seem to have been realized.

{\em Remark} 2. The lattice with $a=b=\sqrt{\pi/\alpha}$ 
is called the 
{\em von Neumann lattice}, since von Neumann 
claimed (without proof) that it is a set of uniqueness 
\cite{ne};
many proofs have later been given [2, 14, 1, 19].
See \cite{grwa} for an attempted repair of the ``defect'' 
of the
von Neumann lattice that it is neither a set of sampling 
nor one
of interpolation.
\vspace{1 mm}

We consider also the analogues in our setting of the 
problems
treated in \cite{be}. We introduce then the 
Banach space $F_{\alpha}^{\infty}$,
consisting of those entire functions $f(z)$ for which
\[
\| f\|_{\infty}=
\| f \|_{\alpha,\infty} 
=\sup_{z}e^{-\alpha|z|^2/2}|f(z)|<\infty.
\]
$\Gamma$ is said to be
a {\em set of sampling} for $F_{\alpha}^{\infty}$ if there 
exists
a positive number $K$ such that
\[
\| f \|_{\infty} \leq K \ 
\sup_{z\in \Gamma}e^{-\alpha|z|^2/2}|f(z)|
\]
for all $f\in F_{\alpha}^{\infty}$. If to every bounded 
sequence
$\{a_{j}\}$ of complex numbers there exists an $f\in
F_{\alpha}^{\infty}$ such that 
$e^{-\alpha|z_{j}|^2/2}f(z_{j})=a_{j}$
for all $j$, we say that $\Gamma=\{z_{j}\}$ is a {\em set 
of interpolation}
for $F_{\alpha}^{\infty}$.
We have then the following counterparts of Beurling's two
density theorems in \cite{be} (we are using
the term sampling instead of balayage as in \cite{be}, 
which seems natural since we no 
longer have the relation to Fourier-Stieltjes transforms). 

\begin{theo}
A discrete set $\Gamma$ is a set of sampling for
$F_{\alpha}^{\infty}$ if and only if it contains a 
uniformly discrete
subset ${\Gamma}'$ for which $D^{-}({\Gamma}')>\alpha/\pi$.
\label{usampling}
\end{theo}

\begin{theo}
A discrete set $\Gamma$ is a set of interpolation for
$F_{\alpha}^{\infty}$ if and only if it is uniformly 
discrete
and $D^{+}({\Gamma})<\alpha/\pi$.
\label{uinterpol}
\end{theo}

Let us  remark, as Beurling did, that the problems and 
some of the
results extend to several variables. We would also like to 
mention the 
following interesting question: What are the corresponding
density theorems
for  weighted Bergman spaces? See \cite{se2,se4} for a 
treatment of this
problem.

\section{The necessity parts of the theorems---indication 
of proof}
In this section we make a few remarks to indicate how to 
prove the 
necessity parts of the theorems. Details are given in 
\cite{se3}.
When unspecified, $p$ is taken to be either $2$ or $\infty$.
 
We remark first that the translations 
\[
(T_{a}f\,)(z)=e^{\alpha \ol{a}z-\alpha|a|^{2}/2}f(z-a)
\]
act isometrically in $F_{\alpha}^{p}$. This translation
invariance implies immediately that 
$\Gamma+z$ is a set of sampling (interpolation) if and 
only if
$\Gamma$ is a set of sampling (interpolation) and it 
permits us to 
translate our analysis around an arbitrary point $z$ to $0$.

Another important feature of $F_{\alpha}^{p}$ is the 
following
compactness property: If $\{f_n \}$ is a sequence in
the ball
\[
\{ f \in F_{\alpha}^{p}:\ \| f\|_p\leq R \},
\]
then there is a subsequence $\{f_{n_k} \}$ converging 
pointwise and
uniformly on compact sets to some
function in the ball. This is immediate from the 
definition of
$F_{\alpha}^{p}$ and a normal family argument.

Following Beurling, for a closed set $\Gamma$, we let 
$W(\Gamma)$ denote
the collection of weak limits of translates $\Gamma+z$
\cite[p. 344]{be}.
The compactness property and the translation invariance of 
$F_{\alpha}^{p}$ make $W(\Gamma)$ a crucial tool in 
our analysis. Indeed, it turns out that all of Beurling's 
arguments
concerning $W(\Gamma)$ can be carried over to our situation.

Most of the work needed to prove the necessity parts of 
Theorems \ref{usampling} and \ref{uinterpol}
consists in transferring Beurling's arguments. In addition 
to
the ingredients mentioned above, a simple substitute for 
Bernstein's 
theorem is used. Moreover, adapting an idea of Landau 
\cite{la}, 
we make use, at a certain stage in the proof of 
Theorem~\ref{usampling}, of the nice properties of the 
normalized monomials  (normalized in $F_{\alpha}^2$),
see [8, p. 39; 15].

For the $L^2$ problem, the basic auxiliary result is the 
following lemma.  

\begin{lem}
There is no discrete subset of ${\bf C}$ that is both a 
set of
sampling and a set of interpolation for $F_{\alpha}^2$.
\label{noriesz}
\end{lem}

This lemma has the following consequences.

\begin{lem}\label{key1}
If $\Gamma$ is a set of sampling for $F_{\alpha}^{2}$, then
so is $\Gamma\setminus\{\zeta\}$ for any $\zeta\in\Gamma$.
\end{lem}

\begin{lem}\label{key2}
If $\Gamma$ is a set of interpolation for $F_{\alpha}^2$,
then so is $\Gamma\cup\{\zeta\}$ for any
$\zeta\not\in\Gamma$.
\end{lem}

The main difficulty in proving the necessity part of
Theorem~\ref{sampling} consists in showing that 
$D^{-}(\Gamma)>\alpha/\pi$
if $\Gamma$ is uniformly discrete and a set of sampling. 
This problem
can now be dealt with in the following way.
Consider such a $\Gamma$. It is easy to
show that $W(\Gamma)$ consists only of sets of
sampling. 
By Lemma~\ref{key1} we have that every set of sampling for
$F_{\alpha}^2$ is a set of uniqueness for 
$F_{\alpha}^\infty$.
For suppose $\Gamma_0$ is a set of sampling for 
$F_{\alpha}^2$ and 
that $g\in F_{\alpha}^{\infty}$ vanishes on
$\Gamma_0$. Then the function
\[
f(z)=g(z)/(z-z_1)(z-z_2),
\]
$z_1,z_2\in \Gamma_0$, belongs to $F_{\alpha}^{2}$ and 
vanishes
on $\Gamma_0\setminus\{z_1,z_2\}$. This contradicts
Lemma~\ref{key1}.

Thus every set in $W(\Gamma)$ is a set of uniqueness for 
$F\,_{\alpha}^{\infty}$. It can be proved	that
$\Gamma$ is a set of sampling for $F_{\alpha}^{\infty}$ if 
and only if 
every $\Gamma_0 \in W(\Gamma)$
is a set of uniqueness for $F_{\alpha}^{\infty}$ (see 
Theorem~3 in 
\cite[p. 345]{be}). Hence by
Theorem \ref{usampling}, $D^{-}(\Gamma)>\alpha/\pi$.

As to the necessity part of Theorem \ref{interpol}, we 
remark that 
Lemma \ref{key2} enables us to carry over Beurling's 
technique used
for the corresponding $L^{\infty}$ problem;  here a slight 
modification of 
the key notion `$\rho(z;\Gamma)$' is needed, see \cite[p. 
352]{be}.

\section{The sufficiency parts of the 
theorems---indication of proof}
Details  can be found in [19].

Let $\Lambda=\{\lambda_{mn}\}$ denote a square lattice; 
that is,
$\lambda_{mn}=\sqrt{\pi/\alpha}\ (m+in)$
for all integers $m,n$ and some positive number $\alpha$. 
$\alpha/\pi$
will be referred to as the {\em density} of $\Lambda$.
We observe that the Weierstrass $\sigma$-function, 
$\sigma(z)$, 
associated to $\Lambda$ plays a role in the 
Bargmann-Fock space analogous to that of the sine in the 
Paley-Wiener
space. This permits us to use techniques similar to some 
of those
employed for functions of exponential type. 

We introduce the analogues of $\sigma(z)$ for uniformly 
discrete
sets that are close to a square lattice in the following 
sense.
$\Gamma=\{z_{mn}\}$ is {\em uniformly close} to $\Lambda$ 
if there exists 
a positive number $Q$ such that 
$|z_{mn}-\lambda_{mn}| \leq Q$
for all $m$ and $n$. 
 
To $\Gamma$, uniformly close to a square lattice, we 
associate a
function $g(z)$ defined by
\begin{equation}
g(z)=(z-z_{00}){\prod_{m,n}}^{\prime}\left(1-%
\frac{z}{z_{mn}}\right)
\exp\left(
\frac{z}{z_{mn}}+\ha\frac{z^{2}}{\lambda^{2}_{mn}}\right)
\label{canon}
\end{equation}
where $z_{00}$ is the point of $\Gamma$ closest to $0$. 
Using the quasi-periodicity of the $\sigma$-function, we 
obtain  the following estimates on the growth of $g$.

\begin{lem}
Let $\Gamma$ be uniformly close to the square lattice 
$\Lambda$ of
density $\alpha/\pi$. Then there exist
constants $C_{1}$, $C_{2}$ and $c$, depending only on $Q$ 
and $q$,
such that for every $z$ we have
\[
|e^{-\alpha|z|^2/2}g(z)|\ \geq \ C_{1}
\ e^{-c |z| \log |z|}\ \mathrm{dist}(z,\Gamma),
\]
\[
|e^{-\alpha|z|^{2}/2} g(z)|\ \leq \ C_{2}\   e^{c |z| \log 
|z|},
\]
and for every $z_{mn}\in\Gamma$ we have
\[
|e^{-\alpha|z_{mn}|^{2}/2}g'(z_{mn})| \ \geq \ C_{1} 
\ e^{-c |z_{mn}| \log |z_{mn}|}.
\]
\label{glem}
\end{lem}

\vskip -12pt
By this lemma and the calculus of residues, we obtain the
following Lagrange-type interpolation formula.

\begin{lem}
Let  $\Gamma=\{z_{mn}\}$ be 
uniformly close to the square lattice of density 
$\beta/\pi$, and let
$g$ be the function
associated to  $\Gamma$ by (\ref{canon}).
If $\alpha<\beta$ we have for each $f\in 
F_{\alpha}^{\infty}$
\[
f(z)=\sum_{m,n}\frac{f(z_{mn})}{g'(z_{mn})}\ 
\frac{g(z)}{z-z_{mn}}
\]
with uniform convergence on compact sets. 
\label{basic}
\end{lem}

The difficulty in proving the sufficiency part of 
Theorem~\ref{sampling}
consists in verifying the left inequality in (\ref{frame}).
We assume then, without loss of generality, that $\Gamma$ 
is uniformly
discrete and uniformly 
close to the
square lattice of density $\beta/\pi$, where
$D^{-}(\Gamma)=\beta/\pi$; see \cite[p. 356]{be} for an 
argument justifying
this claim.
We write
\[
\int_{{\bf C}}\  |f(z)|^{2}\ d\mu_{\alpha}(z)
=\sum_{k,l}\int_{R}\ |(T_{\lambda_{kl}}f)(z)|^{2}\, 
d\mu_{\alpha}(z) 
\]
where $R=\{z=x+iy:\ |x|<\ha\sqrt{1/\alpha},\ 
|y|<\ha\sqrt{1/\alpha}\}$.
In order to estimate the summands on the right, we use 
Lemma \ref{basic}
to write
\[
(T_{\lambda_{kl}}f)(z)=
\sum_{m,n}\frac{(T_{\lambda_{kl}}f)(z_{mn}+\lambda_{kl})}
{g'_{\lambda_{kl}}(z_{mn}+\lambda_{kl})}\ 
\frac{g_{\lambda_{kl}}(z)}{z-z_{mn}-\lambda_{kl}},
\]
where  $g_{\lambda_{kl}}$ is the
function associated to $\Gamma+\lambda_{kl}$ by 
(\ref{canon}). 
Lemma~\ref{glem} is used to estimate this expression, and 
after some 
computation we  obtain the 
desired estimate.

The sufficiency part of Theorem~\ref{usampling} can  be 
proved 
in the same way, or more easily, by Beurling's method 
\cite[p. 346]{be}.

In order to prove the sufficiency parts of Theorems
\ref{interpol} and \ref{uinterpol}, we note first that we 
may 
assume that $\Gamma$
is uniformly close to a square lattice of density 
$\beta/\pi$, where
$D^{+}(\Gamma)=\beta/\pi$ \cite[p. 356]{be}. 
The interpolation problem is then
solved explicitly by the following formula, 
\[
f(z)=\sum_{m,n}a_{mn}e^{\alpha\ol{z}_{mn}z-%
\alpha|z_{mn}|^{2}}
\frac{g_{-z_{mn}}(z-z_{mn})}{(z-z_{mn})}.
\label{intform}
\]
Lemma~\ref{glem} is used to verify  this assertion.

\end{document}